\documentclass[12pt]{amsart}
\usepackage{amssymb}
\usepackage{latexsym}
\usepackage[pdftex]{graphicx}
\usepackage{dsfont}

\newtheorem{theorem}{Theorem}
\newtheorem{lemma}{Lemma}
\newtheorem{proposition}{Proposition}

\newtheorem{corollary}{Corollary}

\newcommand{\beq}{\begin{equation}}
\newcommand{\eeq}{\end{equation}}
\newcommand{\beqna}{\begin{eqnarray*}}
\newcommand{\eeqna}{\end{eqnarray*}}
\newcommand{\beqn}{\begin{equation*}}
\newcommand{\eeqn}{\end{equation*}}
\newcommand{\bp}{\begin{proof}}
\newcommand{\ep}{\end{proof}}
\newcommand{\bprop}{\begin{proposition}}
\newcommand{\eprop}{\end{proposition}}
\newcommand{\bt}{\begin{theorem}}
\newcommand{\et}{\end{theorem}}
\newcommand{\bex}{\begin{Example}}
\newcommand{\eex}{\end{Example}}
\newcommand{\bc}{\begin{corollary}}
\newcommand{\ec}{\end{corollary}}
\newcommand{\bl}{\begin{lemma}}
\newcommand{\el}{\end{lemma}}

\renewcommand{\ge}{\geqslant}
\renewcommand{\le}{\leqslant}

\begin{document}

\title{ Generalized Gr\"unbaum inequality}

 \author{M. ~Meyer, F. ~Nazarov, D.~Ryabogin, and V.~Yaskin}

\address{Universit\'e Paris-Est, Laboratoire d'Analyse et de Math\'ematiques Appliqu\'ees UMR 8050, UPEMLV, UPEC, CNRS
F-77454, Marne-la-Vall\'ee, France }
\email{mathieu.meyer@u-pem.fr}

\address{Fedor Nazarov, Department of Mathematical Sciences, Kent State University, Kent, OH 44242, USA}
 \email{nazarov@math.kent.edu}
 \address{Dmitry Ryabogin, Department of Mathematical Sciences, Kent State University, Kent, OH 44242, USA}
 \email{ryabogin@math.kent.edu}
 \address{Vladyslav Yaskin, Department of Mathematical and Statistical Sciences, University of Alberta, Edmonton, Alberta, T6G 2G1, Canada}
 \email{yaskin@ualberta.ca}

 \subjclass{Primary 52C07, 52A20}

\keywords{log-concave function, Gr\"unbaum inequality}

\thanks{The second and the  third  named authors are
supported in part by U.S.~National Science Foundation Grants  DMS-0800243 and
DMS-1600753, the  fourth author is  partially supported by a grant from NSERC}

 \begin{abstract}
Let $f$ be an integrable  log-concave function on ${\mathbb R^n}$ with the  center of mass at the origin. We show that
$\int\limits_0^{\infty}f(s\theta)ds\ge e^{-n}\int\limits_{-\infty}^{\infty}f(s\theta)ds$ for every $ \theta\in S^{n-1}$, and the constant $e^{-n}$ is the best possible.
\end{abstract}

\maketitle

\section{Introduction}

The classical Gr\"unbaum inequality asserts that for every convex body $K\subset {\mathbb R^d}$ with the center of mass at the origin, and every half-space $H$ whose boundary plane contains the origin, we have
$$
\operatorname{vol}_d(K\cap H)\ge \Big(\frac{d}{d+1}\Big)^d\operatorname{vol}_d(K),
$$
and the equality is attained for any  cone and the half-space containing the vertex of the cone  bounded by the hyperplane parallel to the cone base.

The functional version of the Gr\"unbaum  inequality (due to Lov\'asz and Vempala, \cite{LV}) is
$$
\int\limits_0^{\infty}f(x)dx\ge e^{-1}\int\limits_{-\infty}^{\infty}f(x)dx,
$$
 and the constant $e^{-1}$ is sharp. Here $f$ is any non-negative  integrable log-concave function such that $\int\limits_{-\infty}^{\infty}xf(x)dx=0$.

The following question was asked in \cite{FMY}. Suppose that $L$ is a plane of codimension $n-1$ ($n\ge 2$) passing though the origin and $H_L\subset L$ is a half-plane with the origin on its boundary.
How small can
$\operatorname{vol}_{d-n+1}(K\cap H_L)$ be compared to $\operatorname{vol}_{d-n+1}(K\cap L)$?

This question has a direct analogue  for log-concave functions as well: Find the best constant $a$ such that the inequality
\begin{equation}\label{mda1}
\int\limits_0^{\infty}f(s\theta)ds\ge a\,\int\limits_{-\infty}^{\infty}f(s\theta)ds
\end{equation}
holds for  every non-negative integrable log-concave function $f:{\mathbb R^n}\to{\mathbb R}$ satisfying $\int\limits_{\mathbb R^n}xf(x)dx=0$ and for every  unit vector $\theta$.
In this note we show that $a=e^{-n}$.

The relation between the functional and the body versions is that the best constant in the body version is always not less than the constant in the functional one and they are asymptotically equal when the ambient dimension $d$ tends to infinity.

\section{An example}

Consider the function $f(x)=\mathds{1}_C(x)e^{-x_1}$, where
$$
C=\{x\in{\mathbb R^n}: x_1\ge -n, \, x_2^2+\dots+x_n^2\le(x_1+n)^2\}
$$
is the infinite cone with the vertex at $(-n,0,\dots,0)$ and the axis along the vector $\theta=(1,0\dots,0)$. It is easy to see that $f$ is log-concave and $\int\limits_{\mathbb R^n}xf(x)dx=0$. Indeed, any coordinate except the first one integrates to $0$ by symmetry and
\begin{multline*}
\int\limits_{\mathbb R^n}x_1f(x)dx=c\int\limits_{-n}^{\infty}s(s+n)^{n-1}e^{-s}ds\\
=c\Big( \int\limits_{-n}^{\infty}(s+n)^{n}e^{-s}ds-n \int\limits_{-n}^{\infty}(s+n)^{n-1}e^{-s}ds\Big)\\=ce^n(n!-n(n-1)!)=0.
\end{multline*}

On the other hand, we have
$$
 \int\limits_{-\infty}^{\infty}f(s\theta)ds=\int\limits_{-n}^{\infty}e^{-s}ds=e^n,\qquad \textrm{and}\qquad
 \int\limits_{0}^{\infty}f(s\theta)ds=\int\limits_{0}^{\infty}e^{-s}ds=1.
$$

Thus, the best possible value of $a$ does not exceed $e^{-n}$. It remains
to show that  (\ref{mda1}) always holds with the constant $a=e^{-n}$.  Everywhere below we shall assume that $f$ is positive in the entire space, continuous, strictly logarithmically concave and decays to zero at infinity faster than any exponential function. The general case can be reduced to this one by standard density arguments.

\section{The one-dimensional case}

The case $n=1$ is well-known but we shall still present a proof to motivate the constructions in the following sections. Note that the proof we outline in this section is
by no means the shortest one. However it has the advantage of being ``natural" enough (at least for us) to provide the guidelines for the more involved argument in the multidimensional case.

Let $f:{\mathbb R}\to {\mathbb R}$ be a  positive continuous strictly log-concave   function decaying to zero at infinity faster than any exponent
 and with the center of mass at the origin, i.e.,
$\int\limits_{-\infty}^{\infty}xf(x)dx=0$. Then  we can choose $\beta>0$ such that
$\int\limits_{0}^{\infty}f(0)e^{-\beta x}dx=\int\limits_{0}^{\infty}f(x)dx$. Note that, due to the strict  log-concavity of $f$ and the fast decay property, the graphs of $f$ and $x\mapsto f(0)e^{-\beta x}$ intersect in the way shown on Figure 1;

\begin{figure}[ht]
\includegraphics[height=5.0cm]{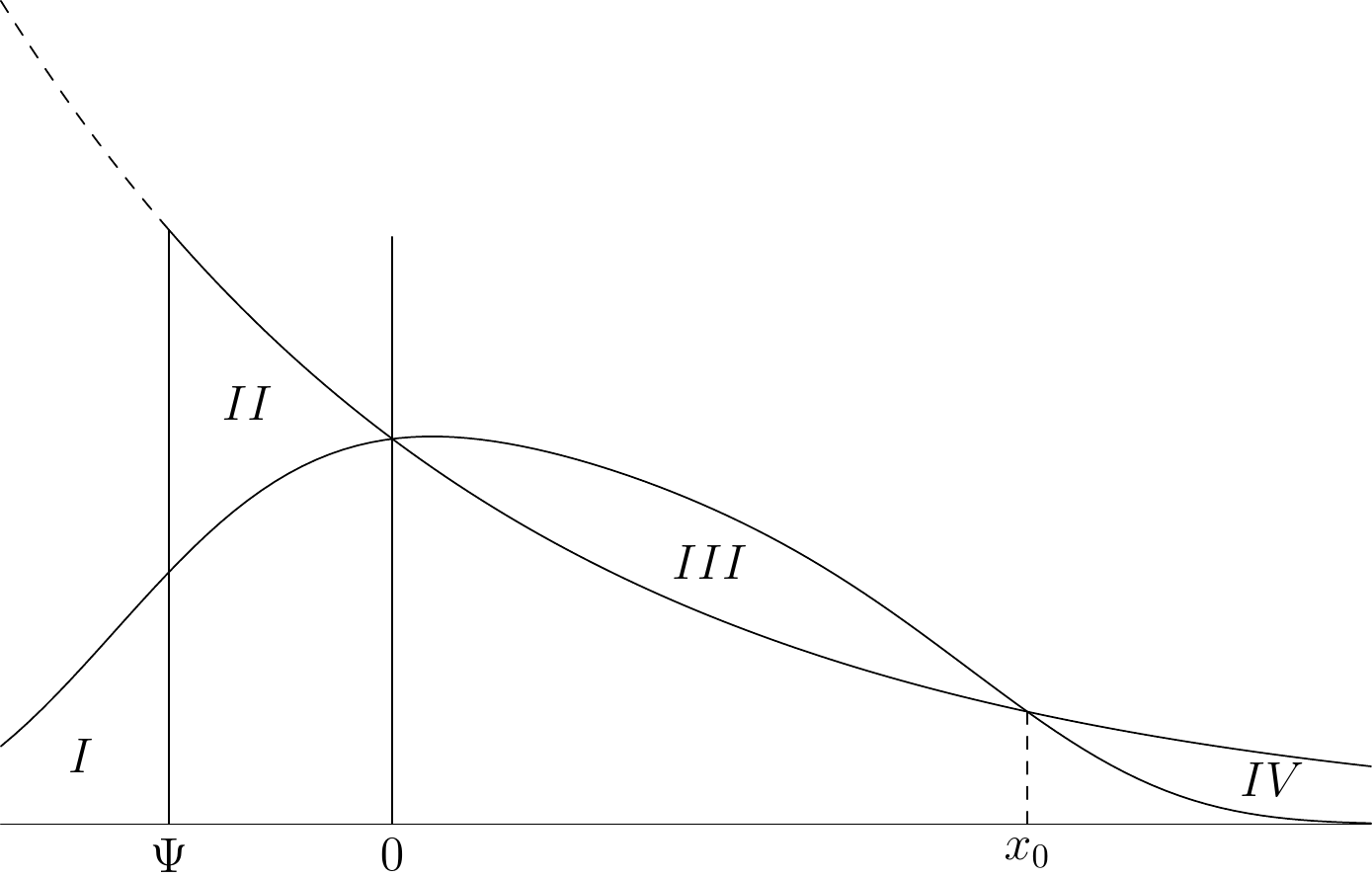}
\caption{The comparison function}
\label{fig1}
\end{figure}

\noindent
i.e., there is $x_0>0$ such that $f(x)\le f(0)e^{-\beta x}$ on $(-\infty,0]\cup[x_0,\infty)$ and $f(x)\ge f(0)e^{-\beta x}$ on $[0,x_0]$.

Define the new function $F(x)=\mathds{1}_{[\Psi,+\infty)}(x)f(0)e^{-\beta x}$, where $\Psi<0$ is chosen so that $\int\limits_{-\infty}^{\infty}F(x)dx=\int\limits_{-\infty}^{\infty}f(x)dx$.
Notice that we also have $\int\limits_{0}^{\infty}F(x)dx=\int\limits_{0}^{\infty}f(x)dx$, and, thereby, $\int\limits_{-\infty}^0F(x)dx=\int\limits_{-\infty}^{0}f(x)dx$. Observe that  the center of mass of $F$ lies to the right of the origin. This is clear from the physical point of view because, when switching from $f$ to $F$ , we move the mass from area $I$ to area $II$ and from area $III$ to area $IV$, i.e., always to the right.

The formal computation is
$$
\int\limits_{-\infty}^{\infty}x(F(x)-f(x))dx=\int\limits_{-\infty}^0x(F(x)-f(x))dx+\int\limits_{0}^{\infty}x(F(x)-f(x))dx=
$$
$$
\int\limits_{-\infty}^0(x-\Psi)(F(x)-f(x))dx+\int\limits_{0}^{\infty}(x-x_0)(F(x)-f(x))dx\ge 0
$$
(both integrands are non-negative).

Thus, if we shift the function $F$ to the left so that the center of mass moves to the origin, we will diminish the integral
$\int\limits_0^{\infty}F(x)dx$ without changing $\int\limits_{-\infty}^{\infty}F(x)dx$, so the constant $a$ we can use for $f$ is at least as large as the one that can be used for $F$. However, $F$ is a pure truncated exponent and it is easy to check that the corresponding ratio of the integrals is $e^{-1}$ regardless of the values of $\beta$ and $f(0)$.

\section{Replacing the slices $f_{x'}(s)$ by exponential functions}

Let us now turn to the multidimensional case. We shall apply the construction of the previous section to the slice $f_0(s)=f(s \theta)$ of the function $f$, i.e., we shall choose $\beta>0$ such that $\int\limits_{0}^{\infty}f(s \theta)ds=\int\limits_{0}^{\infty}f(0)e^{-\beta s}ds$ and get the comparison function $F_0(s)=\mathds{1}_{[\Psi(0),+\infty)}(s)f(0)e^{-\beta s}$. Now we shall fix this value of $\beta$ for the rest of the argument. Our aim will be to make a similar replacement of every slice $f_{x'}(s)=f(x'+s \theta)$, $x'\in \theta^{\perp} \cong{\mathbb R^{n-1}}$
by a function $F_{x'}(s)=\mathds{1}_{[\Psi(x'),+\infty)}(s)H(x')e^{-\beta s}$. Note that we shall use the same $\beta>0$ for every slice. However we are still free to choose any ``height" $H(x')$ we want. Our goal will be to ensure that the resulting function $F(x)=F_{x'}(x_1)$ will still be log-concave and its center of mass will lie on the line $\ell_\theta=\{s\theta:\,s\in {\mathbb R}\}$ to the right of the origin.

That the center of mass of $F$ still lies on the line is automatic because we have not changed the total mass of any slice. To make sure that it lies to the right of the origin, it is enough to ensure that the center of mass of every slice moves to the right after the replacement. That will be achieved by choosing the height function $H$ in an appropriate way.

\section{The critical height}

Let $f:{\mathbb R}\to {\mathbb R}$ be a positive continuous strictly log-concave function tending to zero at infinity faster than any exponential function. The critical height of $f$ is the unique number $H$ such that for $h\ge H$ we have
$\int\limits_{a}^{\infty}f(x)dx\le \int\limits_{a}^{\infty}he^{-\beta x}dx$ for all $a\in{\mathbb R}$, but for every $h<H$ there exists $a\in {\mathbb R}$ such that
$\int\limits_{a}^{\infty}f(x)dx> \int\limits_{a}^{\infty}he^{-\beta x}dx$. In other words, $H=\max\limits_a\beta e^{\beta a}\int\limits_a^{\infty}f(x)dx$.

The critical height is quite easy to visualize geometrically. It is just the number $H$ such that the integrals of $f$ and $He^{- \beta x}$ from the $x$-coordinate of the left intersection point of their graphs to $+\infty$ are equal (see Figure 2).

\begin{figure}[ht]
\includegraphics[height=5.0cm]{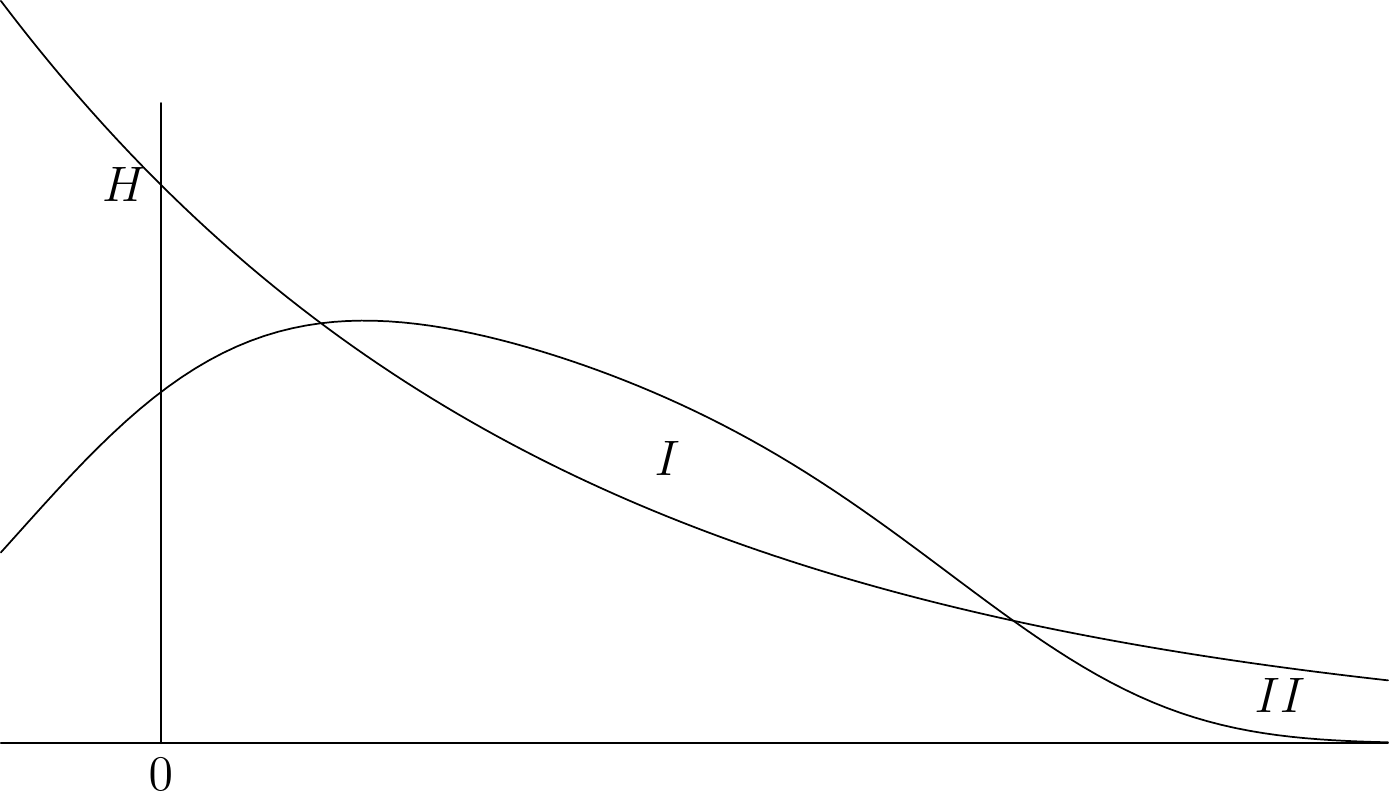}
\caption{The critical height $H$}
\label{fig2}
\end{figure}

The key property of the critical height is that if $h\ge H$ and $\Psi\in {\mathbb R}$ is such that $\int\limits_{\Psi}^{\infty}he^{-\beta x}dx=\int\limits_{-\infty}^{\infty}f(x)dx$, then the center of mass of the function $F(x)=\mathds{1}_{[\Psi,+\infty)}(x)he^{-\beta x}$ is to the right of that of $f$. Again, it is obvious from the physical standpoint because we can move all the excessive mass in  area I to area II. For the formal computations, just use the integration by parts:
$$
\int\limits_{-\infty}^{\infty}x(F(x)-f(x))dx=\int\limits_{-\infty}^{\infty}da\int\limits_{a}^{\infty}(F(x)-f(x))dx\ge 0
$$
because the inner integral is always non-negative by the definition of the critical height.

Now let $H_0(x')$ be the critical height of $f_{x'}$. Then $H_0(0)=f(0)$. We claim that $H_0(x')$ is a continuous logarithmically concave function of $x'$. Indeed, the continuity of $H_0$ follows from the continuity and the fast decay of $f$. Suppose that $x_1'$ and $x_2'$ are any two points in $\theta^{\perp}$ and $h_1<H_0(x_1')$, $h_2<H_0(x_2')$. Then there exist $a_1$, $a_2\in {\mathbb R}$ such that
$$
\int\limits_{a_j}^{\infty}f(x_j'+s\theta)ds>\int\limits_{a_j}^{\infty}h_je^{-\beta s}ds,\qquad j=1,2.
$$
However,  since  $f$ is log-concave, by the Pr\'ekopa-Leindler inequality (see \cite{Ba}, Lecture 5) we have
$$
\int\limits_{\frac{a_1+a_2}{2}}^{+\infty}f\Big(\frac{x_1'+x_2'}{2}+s\theta\Big)ds\ge
\sqrt{\int\limits_{a_1}^{+\infty}f(x_1'+s\theta)ds\int\limits_{a_2}^{+\infty}f(x_2'+s\theta)ds}>
$$
$$
\sqrt{\int\limits_{a_1}^{+\infty}h_1e^{-\beta s}ds\int\limits_{a_2}^{+\infty}h_2e^{-\beta s}ds}=\int\limits_{\frac{a_1+a_2}{2}}^{+\infty}\sqrt{h_1h_2}e^{-\beta s}ds,
$$
which means that $\sqrt{h_1h_2}<H_0(\frac{x_1'+x_2'}{2})$. Since $h_1<H_0(x_1')$ and $h_2<H_0(x_2')$ are arbitrary, we conclude that
$H_0(\frac{x_1'+x_2'}{2})\ge \sqrt{H_0(x_1')H_0(x_2')}$.

Thus, we can  choose a linear function $L:\theta^{\perp}\to{\mathbb R}$ such that $H_0(x')\le f(0)e^{L(x')}$ for every $x'\in\theta^{\perp}$. Now, if we choose the height function $H(x')=f(0)e^{L(x')}$, we will move the center of mass to the right in every slice when replacing $f_{x'}$ by $F_{x'}$.

\section{Convexity and growth of $\Psi$}

Recall that $\Psi(x')$ is defined by the equation
$$
\frac{H(x')}{\beta}e^{-\beta\Psi(x')}=\int\limits_{\Psi(x')}^{\infty}H(x')e^{-\beta s}ds=\int\limits_{-\infty}^{\infty}f(x'+s\theta)ds.
$$
Thus,
$$
\Psi(x')=-\frac{1}{\beta}\Big(\log\frac{\beta}{f(0)}-L(x')+\log\Big(\int\limits_{-\infty}^{\infty}f(x'+s\theta)ds \Big)    \Big).
$$
Since $L$ is linear and $x'\mapsto \int\limits_{-\infty}^{\infty}f(x'+s\theta)ds$ is log-concave, we conclude that $\Psi$ is convex. Moreover, since $f$ decays to zero at infinity faster than any exponent, we see that $\Psi(x')$ grows to $+\infty$ faster than any linear function as $x'\to\infty$. Hence, the comparison function
$F(x)=\mathds{1}_{[\Psi(x'),+\infty)}(x_1)f(0)e^{L(x')-\beta x_1}$ is supported by the convex domain $\{x\in{\mathbb R^n}:x_1\ge \Psi(x')\}$ and is purely exponential in that domain. Moreover, the center of mass of $F$ is on the line $\ell_{\theta}$ to the right of the origin and
$\int\limits_{-\infty}^{\infty}F(s\theta)ds=\int\limits_{-\infty}^{\infty}f(s\theta)ds$, $\int\limits_{0}^{\infty}F(s\theta)ds=\int\limits_{0}^{\infty}f(s\theta)ds$. Thus, if we move $F$ in the direction $-\theta$ so that the center of mass moves to the origin, then the ratio of the corresponding integrals for $F$ will be less than or equal to the one for $f$.
Therefore, it is enough to investigate the functions of the type
$F(x)=\mathds{1}_{Q}(x)f(0)e^{L(x')-\beta x_1}$, where $Q=\{x\in{\mathbb R^n}:x_1\ge \Psi(x')\}$ is an unbounded convex domain  lying to the right of the graph of a convex function $\Psi:\theta^{\perp}\to {\mathbb R}$ growing faster than any linear function at infinity and such that $F(x)$ has the center of mass at the origin.

\section{Symmetrization}

Let us now consider the ``level planes" $\{x\in{\mathbb R^n}:L(x')-\beta x_1=\operatorname{const}\}$. Since $\Psi$ grows faster than any linear function, all sections of $Q$ by the level planes are bounded.
We shall symmetrize $Q$ by replacing any such section by the $(n-1)$-dimensional ball lying in the same plane and of the same $(n-1)$-dimensional volume  centered on the line $\ell_\theta $ (see Figure 3).

\begin{figure}[ht]
\includegraphics[height=5.0cm]{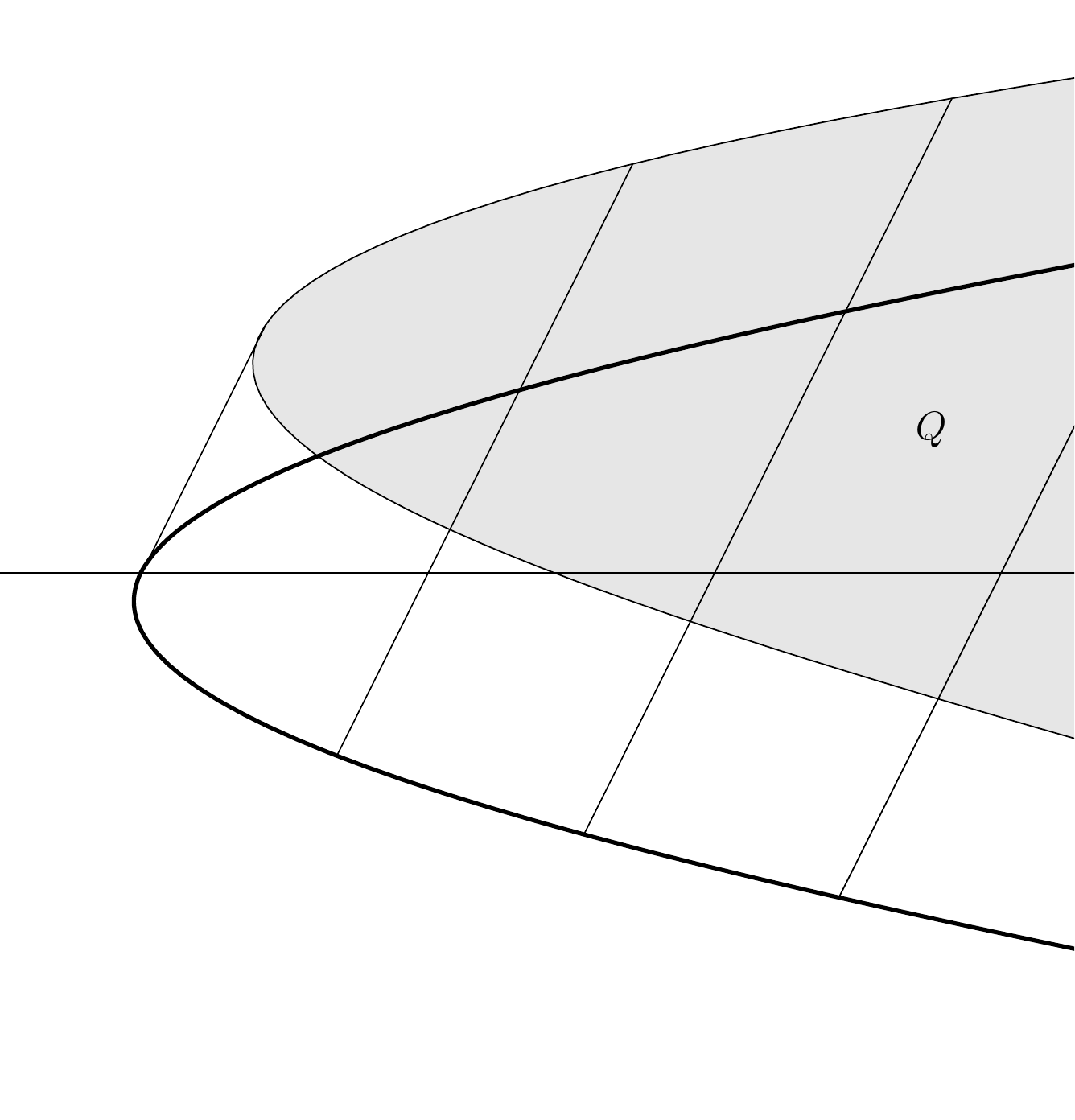}
\caption{Symmetrization}
\label{fig3}
\end{figure}

Note that  this operation leaves the center of mass at the origin. Indeed, since every section  of the symmetrized body by a level plane has its center of mass on $\ell_\theta$, the center of mass of the whole body  stays on $\ell_\theta$. On the other hand, since the masses of sections of $Q$  by  level planes are preserved, it can move only in the direction of the level planes, which is transversal to $\ell_\theta$.

Also, after this operation the integral
$\int\limits_{0}^{\infty}F(s\theta)ds$ stays the same while the integral $\int\limits_{-\infty}^{\infty}F(s\theta)ds$ may only increase (due to the extension of the support of the restriction of $F$ to the line $\ell_\theta$). We can also apply an appropiate linear transformation that makes the level planes perpendicular to $\theta$. So we need only to consider the functions $F(x)=\mathds{1}_{Q}(x)e^{-\beta x_1}$, where $Q$ is an unbounded (from the right) convex body of revolution around the axis $\ell_\theta $.

\section{The comparison cone}

Let now $s_0\theta$ ($s_0<0$) be the point on the boundary of the body of revolution $Q$. Let $B$ be the cross-section of $Q$ by the plane $\theta^{\perp}$ (see Figure 4).

\begin{figure}[ht]
\includegraphics[height=7.0cm]{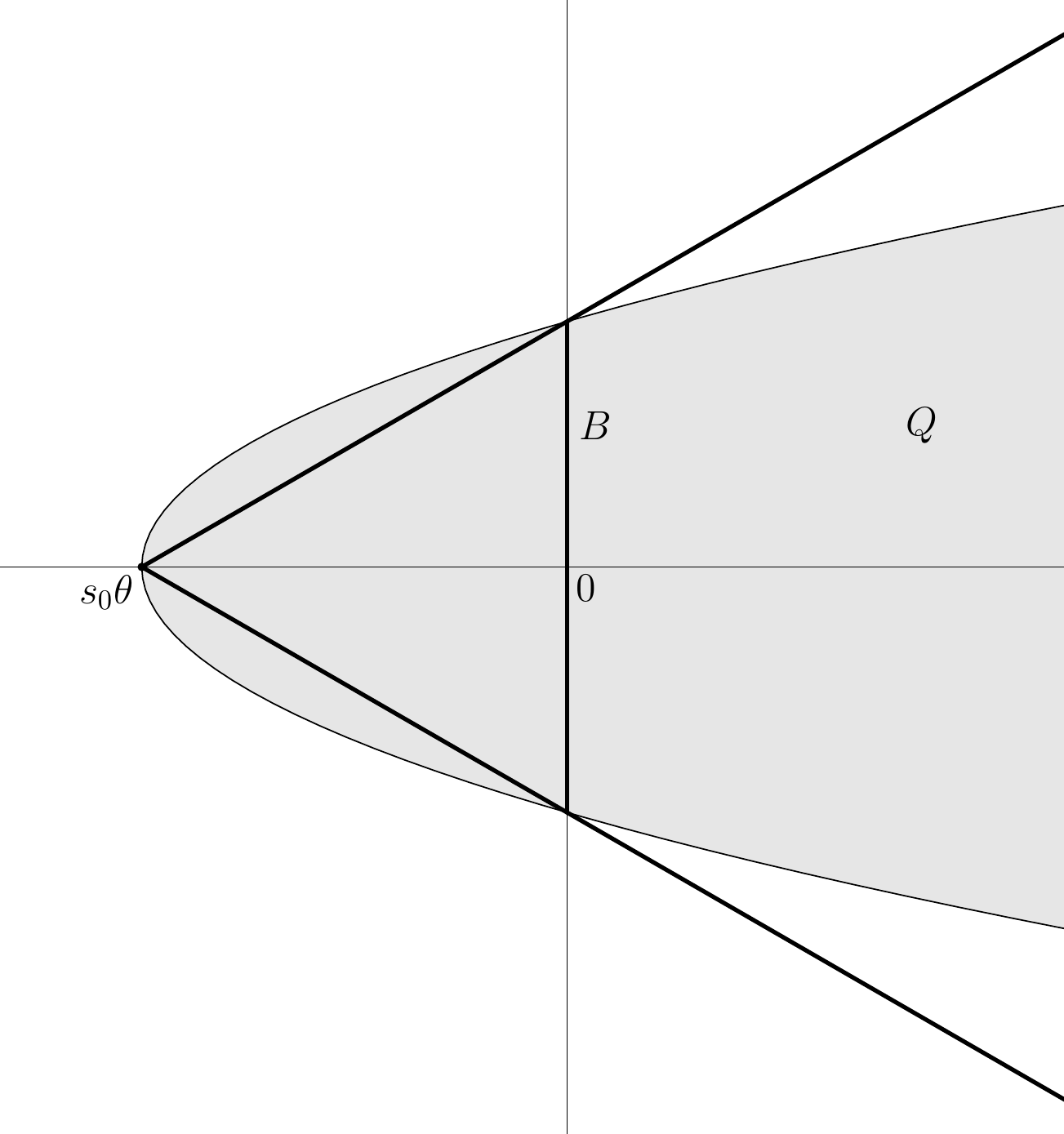}
\caption{The comparison cone}
\label{fig4}
\end{figure}

Consider the infinite cone $C$ with the vertex $s_0\theta$ whose cross-section by $\theta^{\perp}$ is also $B$ and the function
$\widetilde{F}(x)=\mathds{1}_{C}(x)e^{-\beta x_1}$. The restriction of $\widetilde{F}$ to $\ell_\theta$ coincides with that of $F$, so to finish the proof, it will suffice to show that the center of mass of $\widetilde{F}$ lies to the right of the origin. But we have $x_1(\widetilde{F}(x)-F(x))\ge 0$ and, thereby,
$$
\int\limits_{{\mathbb R^n}}x_1\widetilde{F}(x)dx=\int\limits_{{\mathbb R^n}}x_1(\widetilde{F}(x)-F(x))dx\ge 0.
$$

\end{document}